\documentclass[11pt,reqno]{article}
\usepackage[all]{xy}
\usepackage{amsmath}
\usepackage{amsthm}
\usepackage{amscd}
\usepackage{amssymb}

\input{epsf.tex}

\numberwithin{equation}{section}

\newcommand{\nb}[1]{#1\nobreakdash-}

\theoremstyle{definition}

\theoremstyle{plain}
\newtheorem{theorem}{Theorem}
\newtheorem{proposition}[theorem]{Proposition}
\newtheorem{lemma}[theorem]{Lemma}
\newtheorem{corollary}[theorem]{Corollary}
\newtheorem{claim}[theorem]{Claim}

\newcounter{remarks}
{\paragraph*{Remarks}\smallskip
     \begin{list}{\arabic{remarks}. }{\usecounter{remarks}%
          \setlength{\leftmargin}{0in}%
          \setlength{\rightmargin}{0in}%
          \setlength{\labelsep}{0pt}%
          \setlength{\labelwidth}{0pt}%
          \setlength{\listparindent}{0pt}%
     }
}
{
\end{list}
}

\newcommand\ds\displaystyle
\newcommand\wt[1]{\widetilde{#1}}

\DeclareMathOperator{\Out}{Out}

\DeclareMathOperator{\Fix}{Fix}

\DeclareMathOperator{\Homeo}{Homeo}
\DeclareMathOperator{\QIMap}{\widehat{\QI}}

\DeclareMathOperator{\Isom}{Isom}

\DeclareMathOperator{\Aut}{Aut}
\DeclareMathOperator{\Stab}{Stab}

\newcommand\R{{\mathbf R}}
\newcommand\reals{\R}

\renewcommand\H{{\mathbf H}}
\newcommand\Hyp\H
\newcommand\hyp{\H}
\newcommand\Z{{\mathbf Z}}

\newcommand\solv{{\ensuremath{\text{\scshape solv}}}}
\newcommand\Solv\solv

\newcommand\inject{\hookrightarrow}

\newcommand\infinity{\infty}
\newcommand\bndry{\partial}
\newcommand{\bdy}{\bndry}
\newcommand{\from}{\colon}
\def\composed{\circ}
\newcommand\suchthat{\bigm|}
\newcommand\inverse{{-1}}
\newcommand\inv{\inverse}

\newcommand\union{\cup}

\newcommand\norm[1]{\left\|  #1 \right\|}

\newcommand\Id{\text{Id}}

\newcommand\intersect{\cap}

\newcommand\subgroup{<}
\newcommand\semidirect{\rtimes}

\newcommand\Teichmuller{Teichm\"uller}

\newcommand\MCG{\mathcal{MCG}}

\DeclareMathOperator\QI{QI}
\newcommand\cross{\times}

\newcommand\Haus{{\mathcal H}}

\renewcommand\O{{\mathcal O}}

\renewcommand\P{{\mathbf P}}

\newcommand\fol{{\mathcal F}}

\newcommand\F\fol

\newcommand\C{\mathcal C}
\newcommand\Teich{\mathcal T}
\newcommand\T\Teich

\newcommand\Porb{\mathcal P}
\newcommand\Qorb{\mathcal Q}

\newcommand\PMF{{\P\mathcal MF}}
\newcommand\MF{{\mathcal MF}}

\DeclareMathOperator\QSym{QSym}

\newcommand\e\epsilon
\newcommand\<\langle
\renewcommand\>\rangle
\DeclareMathOperator\VN{VN}
\DeclareMathOperator\VC{VC}

\newcommand\splitting\succ
\newcommand\splitsto\splitting

\title{Fiber respecting quasi-isometries of \\ surface group extensions}
\author{Lee Mosher\thanks{Supported in part by NSF grant
DMS 0103208.}}

\begin{document}
\maketitle

\begin{abstract}
Let $S$ be a closed, oriented surface of genus $g \ge 2$, and consider the
extension $1 \to \pi_1 S \to \MCG(S,p) \to \MCG(S) \to 1$, where $\MCG(S)$
is the mapping class group of $S$, and $\MCG(S,p)$ is the mapping class
group of $S$ punctured at $p$. We prove that any quasi-isometry of
$\MCG(S,p)$ which coarsely respects the cosets of the normal subgroup
$\pi_1 S$ is a bounded distance from the left action of some element of
$\MCG(S,p)$. Combined with recent work of Kevin Whyte this implies that if
$K$ is a finitely generated group quasi-isometric to $\MCG(S,p)$ then there
is a homomorphism $K \to \MCG(S,p)$ with finite kernel and finite index
image.  Our work applies as well to extensions of the form $1 \to \pi_1 S
\to\Gamma_H \to H \to 1$, where $H$ is an irreducible subgroup of
$\MCG(S)$---we give an algebraic characterization of quasi-isometries of
$\Gamma_H$ that coarsely respect cosets of $\pi_1 S$.
\end{abstract}

\section{Introduction}

A \emph{surface group extension} $1\to \pi_1(S) \to \Gamma \to H \to 1$ is
a short exact sequence where $S$ is a closed, oriented surface of genus $\ge
2$. When $S$ is fixed the universal such extension has two isomorphic
descriptions, according to the Dehn-Nielsen-Epstein-Baer theorem:
$$\xymatrix{
1 \ar[r] & \pi_1(S) \ar[r] \ar@{=}[d] & \MCG(S,p) \ar[r] \ar@2{~}[d] &
\MCG(S) \ar[r] \ar@2{~}[d] & 1 \\ 
1 \ar[r] & \pi_1(S) \ar[r] & \Aut(\pi_1(S)) \ar[r] & \Out(\pi_1(S)) \ar[r]
& 1 }
$$
The universality property says that every surface group extension
arises from a homomorphism of short exact sequences
$$\xymatrix{
1 \ar[r] & \pi_1(S) \ar[r] \ar@{=}[d] & \Gamma \ar[r] \ar[d] & H
\ar[r] \ar[d] & 1 \\ 
1 \ar[r] & \pi_1(S) \ar[r] & \MCG(S,p) \ar[r] & \MCG(S) \ar[r]
& 1 }
$$
and the extension is uniquely determined, up the the appropriate
equivalence relation, by the homomorphism $H \to \MCG(S)$. In the special
case of an extension determined by inclusion of a subgroup $H \subgroup
\MCG(S)$, the extension group will be denoted $\Gamma_H$.

In studying the large scale geometry of a finitely generated group
$\Gamma$, it is natural to consider the quasi-isometry group $\QI(\Gamma)$, which is
the group of self-quasi-isometries of $\Gamma$ with its word metric, modulo the
relation that two quasi-isometries from $\Gamma$ to itself are equivalent if
their sup distance is finite. When $\Gamma$ is a surface group extension, in
many situations it is also natural to focus on the subgroup $\QI_f(\Gamma)
\subgroup\QI(\Gamma)$ consisting of classes of self quasi-isometries that
coarsely respect the decomposition of $\Gamma$ into cosets of $\pi_1(S)$; we
call these \emph{fiber respecting} quasi-isometries. For example, when $H$
is a free group and so $\Gamma_H$ is surface-by-free, then every
quasi-isometry of $\Gamma$ is fiber respecting, that is, $\QI_f(\Gamma) =
\QI(\Gamma)$ \cite{FarbMosher:ABC}. For another example, the left action of
$\Gamma$ on itself is fiber respecting, by normality of $\pi_1(S)$, and so
the induced homomorphism $\Gamma \to \QI(\Gamma)$ factors through a
homomorphism $\Gamma \to \QI_f(\Gamma)$.

Our goal in this paper is to study the group $\QI_f(\Gamma)$ in various
cases. For example, when $\Gamma$ is the universal extension group
$\MCG(S,p)$ itself we obtain:

\begin{theorem}
\label{TheoremMCG}
The injection $\MCG(S,p) \to \QI_f(\MCG(S,p))$ is an isomorphism. 
\end{theorem}

Kevin Whyte has developed new methods, using uniformly finite homology, for
showing that every self quasi-isometry is fiber respecting. Whyte's methods
apply, for example, to a surface-by-free group, but they also apply in more
general situations situations where the quotient group $H$ has dimension
greater than~1. In particular, Whyte proves:

\begin{theorem}[Whyte] 
\label{TheoremWhyte}
Every self quasi-isometry of $\MCG(S,p)$ is fiber respecting, that is,
$\QI_f(\MCG(S,p))=\QI(\MCG(S,p))$.
\end{theorem} 

As mentioned, Whyte's techniques apply to far more general situations, but
for a proof of just this theorem see \cite{Mosher:Durham03}.

Combining Theorems~\ref{TheoremMCG} and~\ref{TheoremWhyte} we obtain a
strong quasi-isometric rigidity theorem for once-punctured mapping class
groups:

\begin{corollary}[Mosher--Whyte]
\label{CorollaryQIGroup}
The injection $\MCG(S,p) \to \QI(\MCG(S,p))$ is an isomorphism. 
\end{corollary}

Each of these three results has a more quantitative version---in the last
corollary, the conclusion is that for each $K \ge 1$, $C \ge 0$ there
exists $A \ge 0$ such that if $\phi \from \MCG(S,p) \to \MCG(S,p)$ is a
$K,C$ quasi-isometry, then there exists $g \in \MCG(S,p)$ such that
$d(\phi(f),gf) \le A$ for all $f \in \MCG(S,p)$. From this quantitative
version, it follows by a standard argument that:

\begin{corollary}[Mosher--Whyte] 
\label{CorollaryQIRigidity}
If $G$ is any finitely generated group
quasi-isometric to $\MCG(S,p)$ then there is a homomorphism $G \to
\MCG(S,p)$ with finite kernel and finite index image.
\end{corollary}

The kind of quasi-isometric rigidity given in
Corollaries~\ref{CorollaryQIGroup} and~\ref{CorollaryQIRigidity} is
generally conjectured to occur for mapping class groups of nonexceptional
finite type surfaces, and this is the first case where the conjecture is
verified. 

Theorem~\ref{TheoremMCG} will follow from a more general result which
applies to surface group extensions $1 \to \pi_1(S)\to \Gamma_H \to H \to
1$ where $H$ is a subgroup of $\MCG(S)$ that contains at least one
pseudo-Anosov mapping class. For the reader who wants a direct proof of Theorem~\ref{TheoremMCG} without
the generalization to groups $\Gamma_H$, and who wants to see details about
Theorem~\ref{TheoremWhyte}, we refer to \cite{Mosher:Durham03}.

By Ivanov's theorem, a subgroup $H \subgroup \MCG(S)$ contains a
pseudo-Anosov mapping class if and only if $H$ is infinite and irreducible,
meaning that there does not exist an essential curve system on $S$ that is
preserved by every element of $H$.

Consider the extension $\Gamma_H$ of $\pi_1(S)$ by an infinite,
irreducible subgroup $H \subgroup \MCG(S)$. The homomorphism $\Gamma_H \to
\QI_f(\Gamma_H)$ is an injection, but in general there may be additional
elements of $\QI_f(\Gamma_H)$ which can be thought of as ``hidden''
symmetries of $\Gamma_H$. To describe these, let $S \to\O_H$ be the orbifold
covering map of largest degree such that $H$ descends to a subgroup of
$\MCG(\O_H)$, denoted $H'$. Let $\C$ be the relative commensurator of
$H'$ in $\MCG(\O_H)$, that is, the group of all $g\in\MCG(\O_H)$ for which
$g^\inv H' g \intersect H'$ has finite index in both $g^\inv H' g$ and $H'$.
We obtain an extension group
$$1 \to \pi_1(\O_H) \to \Gamma_\C \to \C \to 1
$$
From this description we obtain an obvious injection $\Gamma_\C \to
\QI_f(\Gamma_H)$. The question arises whether there is anything else in
$\QI_f(\Gamma_H)$, and the answer is no:

\begin{theorem} 
\label{TheoremMain}
If $1 \to \pi_1(S) \to \Gamma_H \to H \to 1$ is a surface
group extension where $H$ is an infinite, irreducible subgroup of
$\MCG(S)$, then the injection $\Gamma_\C \to \QI_f(\Gamma_H)$ is an
isomorphism.
\end{theorem}

To derive Theorem~\ref{TheoremMCG} from Theorem~\ref{TheoremMain}, we
simply note that $S$ itself is the maximal subcover of $S$ to which
$\MCG(S)$ descends, and $\MCG(S)$ is its own relative commensurator in
$\MCG(S)$. We will also see how to prove the more quantitative versions of
Theorem~\ref{TheoremMCG} and Corollary~\ref{CorollaryQIGroup}.

In the special case when $H$ is a Schottky subgroup of $\MCG(S)$,
Theorem~\ref{TheoremMain} was proved in \cite{FarbMosher:sbf}. The
general proof of Theorem~\ref{TheoremMain} follows the same outline, but
with several changes in detail, the most important of which we highlight
here. First, whereas the nontrivial elements of a Schottky subgroup are
entirely pseudo-Anosov, the most we can say for an irreducible subgroup is
that it is generated by its pseudo-Anosov elements (with a simple
exception; see Lemma~\ref{PropPsAGenerators}); this distinction permeates
the proof. Second, the definition of ``the orbifold subcover of largest
degree to which a subgroup descends'' was handled more easily in
\cite{FarbMosher:sbf} because we were working only with free subgroups, but
it turns out that descent is sufficiently well behaved even for arbitrary
subgroups (see Lemma~\ref{LemmaExtension}). Third, in \cite{FarbMosher:sbf}
the key step of proving that $\Gamma_\C\to\QI_f(\Gamma_H)$ is surjective
depended on an explicit computation of the relative commensurator of a
Schottky subgroup of a mapping class group; in the present situation the
computation is far less explicit, and so we must work with some general
properties of relative commensurators.

\section{Preliminaries}

\paragraph{Coarse language.} A map $f \from X \to Y$ between two metric
spaces is a \emph{$K,C$ quasi-isometry}, $K \ge 1$, $C \ge 0$, if 
$$\frac{1}{K} d_X(a,b) - C \le d_Y(fa,fb) \le K d_X(a,b) + C
$$
for all $a,b \in X$, and for all $c \in Y$ there exists $a \in X$ such that
$d_Y(fa,c) \le C$. A map $\bar f \from Y \to X$ is a \emph{coarse inverse}
for $f$ if $d_X(\bar f \composed f x,x)$ and $d_Y(f \composed \bar f y,y)$
are uniformly bounded for $x \in X$, $y \in Y$.

Given a metric space $X$, let $\QIMap(X)$ denote the collection of all
self-quasi-isometries of $X$ with the operation of composition. Two
elements $f,g \in \QIMap(X)$ are considered to be equivalent if
$\sup\{d(fx,gx) \suchthat x \in X\} < \infinity$. Composition descends to a
group operation on the quotient $\QI(X)$, whose identity element is the
class of the identity. The inverse of the class of $f \in \QIMap(X)$ is the
class of any coarse inverse $\bar f$.

Given a metric space $X$ and two subset $A,B$, we say that $A$ is
\emph{coarsely included in $B$} if there exists $r \ge 0$ such that $A
\subset N_r(B)$. We say that $A,B$ are \emph{coarsely equivalent} if each
is coarsely contained in the other; equivalently, the Hausdorff distance
$d_\Haus(A,B) = \inf \{r \suchthat A \subset N_r(B), B
\subset N_r(A)\}$ is finite.

We shall need some facts about coarse inclusions among subgroups of a
finitely generated group, taken from \cite{MSW:QTTwo}. Recall that two
subgroups $A,B \subgroup G$ are \emph{commensurable} if $A\intersect B$ has
finite index in $A$ and in $B$.

\begin{proposition}
\label{PropSubgroups}
If $G$ is a finitely generated group with the word metric and $A,B$ are
subgroups of $G$, then $A$ is coarsely contained in $B$ if and only if $A
\intersect B$ has finite index in $A$, and $A,B$ are coarsely equivalent if
and only if $A,B$ are commensurable in $G$. \qed
\end{proposition}

Consider a finitely generated group $G$ and a normal subgroup $N \subgroup
G$. For each $g \in G$, since $g N g^\inv = N$ it follows that the
Hausdorff distance from each coset $gN=Ng$ to the subgroup $N$ is finite,
equal to at most the word length of $g$. A quasi-isometry $\Phi \from G \to
G$ is \emph{fiber preserving} with respect to $N$ if there is a constant $R$
such that for each coset $gN$ there exists a coset $g'N$ such that
$d_\Haus(\Phi(gN),g'H) \le R$. A composition of fiber preserving
quasi-isometries is fiber preserving, and the identity is fiber preserving,
so the subset of $\QI(G)$ represented by fiber preserving quasi-isometries
(with respect to $N$) is a subgroup denoted $\QI_f(G)$. Each fiber
preserving quasi-isometry $\Phi \from G \to G$ induces a quasi-isometry
$\phi \from G/N \to G/N$ which is well-defined up to equivalence: given the
constant $R \ge 0$ as above, for each coset $gN$ choose $g'N$ as above and
define $\phi(gN)=g'N$.

\paragraph{The virtual normalizer and virtual centralizer of a
pseudo-Anosov mapping class.} Consider a pseudo-Anosov $f \in \MCG(S)$ with
stable and unstable foliations $\fol^s$, $\fol^u$. Recall that the action
of $f$ on $\PMF(S)$ has north-south dynamics with repelling fixed point
$\P\fol^s$ and attracting fixed point $\P\fol^u$. The virtual normalizer
$\VN(f)$ and the virtual centralizer $\VC(f)$ of the infinite cyclic
subgroup $\<f\>$ generated by $f$ are defined to be
\begin{align*}
\VN(f) &= \{ g \in \MCG(S) \suchthat g^\inv \<f\> g \intersect \<f\>
\quad\text{has finite index in $g^\inv \<f\> g$ and in $\<f\>$}\} \\
\VC(f) &= \{ g \in \MCG(S) \suchthat \exists n \ge 0 \quad\text{such that}
\quad g^\inv f^n g = f^n \}
\end{align*}
We need the following alternative description of these subgroups, which is
due to McCarthy and can be found in his preprint 
\cite{McCarthy:normalizers}.

\begin{proposition} The subgroup $\VN(f)$ is the stabilizer of the set
$\Fix(f) = \{\P\fol^s,\P\fol^u\}$, and $\VC(f)$ is the kernel of the action
of $\VN(f)$ on this set, and so $\VC(f)$ is a subgroup of index at most~2
in $\VN(f)$. Each of these subgroups is virtually infinite cyclic, in fact,
$\VN(f)$ is the maximal virtually cyclic subgroup of $\MCG(S)$ containing
$f$, and $\VC(f)$ is the maximal subgroup containing $f$ which is virtually
cyclic and contains no $D_\infinity$ subgroup. 
\end{proposition}

One can easily show that the inclusion $\VC(f) \subgroup \VN(f)$ is proper
if and only if there is a pseudo-Anosov element of $\VN(f)$ that is
conjugate to its own inverse.

\begin{proof} For completeness, here is a proof of these facts.

Consider the \Teichmuller\ geodesic $\gamma$ with endpoints
$\P\fol^s,\P\fol^u$, and clearly $\Stab(\gamma)=\Stab(\Fix(f))$. The
geodesic $\gamma$ is the axis of $f$ in \Teichmuller\ space. The group
$\Stab(\gamma)$ acts properly and cocompactly on $\gamma$ and so is a
virtually cyclic group, containing $\<f\>$ with finite index, and therefore
implying that $\Stab(\gamma)\subset \VN(f)$. Conversely, if $g
\in \MCG(S) - \Stab(\Fix(f))$ then $g^\inv(\gamma) \ne \gamma$, and $g^\inv
f^n g$ is a pseudo-Anosov mapping class with axis $g^\inv(\gamma)$ for any
$n$. By uniqueness of axes, for any $m,n \ne 0$ we have $g^\inv f^n g \ne
f^m$. This proves that $\VN(f) = \Stab(\Fix(f))$, and the remaining
contentions are easily obtained.
\end{proof}

\paragraph{Irreducible subgroups of $\MCG(S)$.} Ivanov's theorem gives a
dichotomy for infinite subgroups $G \subgroup \MCG(S)$: either $G$ is
\emph{reducible} meaning that there exists a collection $\C$ of pairwise
disjoint, essential simple closed curves in $S$ such that each element of
$G$ preserves $\C$ up to isotopy; or $G$ contains a pseudo-Anosov element.
Briefly, every infinite, irreducible subgroup of $\MCG(S)$ contains a
pseudo-Anosov element. 

We wish to show that irreducible subgroups are generated by their
pseudo-Anosov elements. Actually, this is not quite true: if the
pseudo-Anosov mapping class $f$ has the property that the inclusion $\VC(f)
\subgroup \VN(f)$ is proper, then the subgroup of $\VN(f)$ generated by the
pseudo-Anosov elements is precisely $\VC(f)$. But this is essentially the
only counterexample, as the following result shows:

\begin{proposition}
\label{PropPsAGenerators}
Consider a subgroup $H \subgroup \MCG(S)$ which contains a pseudo-Anosov
element. 
\begin{description}
\item[The nonelementary case] If $H\not\subgroup\VN(f)$ for any
pseudo-Anosov $f\in\MCG(S)$ then $H$ is generated by its pseudo-Anosov
elements.  
\item[The elementary case] If $H\subgroup\VN(f)$ for some pseudo-Anosov $f
\in \MCG(S)$ then the subgroup of $H$ generated by the pseudo-Anosov
elements is $H \intersect \VC(f)$.
\end{description}
\end{proposition}

\begin{proof} The proof follows quickly from Exercise 3a on page 102 of
\cite{Ivanov:subgroups}, which says that if $f \in \MCG(S)$ is
pseudo-Anosov and if $g \in \MCG(S)$ has the property that $g(\fol^u(f))
\ne \fol^s(f)$, then $f^n g$ is pseudo-Anosov for all sufficiently large $n
\ge 0$. Here are the details, including the solution of the exercise.

Consider first the elementary case, $H \subgroup \VN(f)$. Note that each
pseudo-Anosov element of $\VN(f)$ fixes $\P\fol^s$ and $\P\fol^u$, and so
the subgroup they generate is contained in $\VC(f)$. Conversely, given $g
\in \VC(f) \intersect H$ and a pseudo-Anosov element $f \in \VN(f)
\intersect H$, clearly $h = f^n g$ is pseudo-Anosov for sufficiently large
$n$, and $g = h f^{-n}$.

Consider next the nonelementary case. For any nonidentity element $g \in
\MCG(S)$ we may choose a pseudo-Anosov $f \in \MCG(S)$ so that $g \not\in
\VN(f)$. It follows that $g$ does not stabilize the set
$\{\P\fol^u,\P\fol^s\}$ and so, replacing $f$ by $f^\inv$ if necessary, we
may assume that $g(\fol^u(f)) \ne \fol^s(f)$.

We shall prove that $f^n g$ is pseudo-Anosov for sufficiently large $n$ by
showing that $f^n g$ has no fixed points in $\MF$.

Choose neighborhoods $V^s$ of $\P\fol^s$ and $V^u$ of $\P\fol^u$ so that
$V^s \intersect V^u = \emptyset$, and so that $g^\inv V^s \intersect V^u =
\emptyset$; this is possible because $g(\P\fol^u) \ne \P\fol^s$. Let $W =
g^\inv V^s$. If $n$ is sufficiently large it follows that $f^n g (\PMF -
W) \subset V^u$, and so the only possible fixed points of $f^n g$ in
$\PMF$ are in $W$ or in $V^u$. 

For a subset $A \subset \PMF$ let $\R_+ A$ denote the preimage of $A$ under
the projectivation $\MF \to \PMF$. The only possible fixed points of $f^n
g$ in $\MF$ are in $\R_+ W$ or in $\R_+ V^u$.

Consider the action of $f^n g$ on $\MF$. Choose a continuous norm
$\norm{\cdot}$ on $\MF$, meaning that $\norm{r \cdot \fol} = r \norm{\fol}$
for each $\fol \in \MF$, $r > 0$. Since $f(\fol^u) = \lambda \fol^u$ with
$\lambda>1$, it follows that if $n$ is sufficiently large then there is a
constant $m>1$ such that $\norm{f^n g(\fol)} \ge m \norm{\fol}$ for all
$\fol \in \MF - \R_+ W$; it follows that $f^n g$ has no fixed points in
$\R_+ V^u$. Similarly, if $n$ is sufficiently large then $\norm{g^\inv
f^{-n}(\fol)} \ge m \norm{\fol}$ for all $\fol \in \R_+ W$. It follows that
$f$ has no fixed points in $\R_+ W$.
\end{proof}

\paragraph{\Teichmuller\ space and its canonical bundle.}

Let $\Teich = \Teich(S)$ denote the \Teichmuller\ space of $S$. Let
$X_\Teich \to \Teich$ denote the canonical hyperbolic plane bundle over
$\Teich$ on which $\pi_1 S$ acts: the fiber $D_x$ over a point $x \in
\Teich$ is the universal cover of a hyperbolic surface representing the
point $x$. The action of $\MCG(S)$ on $\Teich(S)$ lifts  to an action of
$\MCG(S,p)$ on $X_\Teich$, indeed $X_\Teich$ can be identified naturally as
the \Teichmuller\ space of the once punctured surface $S-p$
\cite{Bers:FiberSpaces}.

We fix once and for all an $\MCG(S,p)$ equivariant Riemannian metric on
$X_\Teich$ whose restriction to each fiber is the hyperbolic metric on that
fiber. 

A \emph{path} in $\Teich$ will always mean a piecewise geodesic map of a
subinterval of $\reals$ into $\Teich$, and a \emph{bi-infinite path} has
domain $\reals$. By pulling back the fiber bundle $X_\Teich\to\Teich$ to the
domain of any path $\gamma$ we obtain the canonical $\hyp^2$ bundle
$X_\gamma$ over~$\gamma$.

Consider two bi-infinite paths $\gamma,\gamma' \from \reals \to
\Teich$. We say that $\gamma,\gamma'$ are \emph{fellow travellers} if there
exists a quasi-isometric homeomorphism $h \from \reals \to \reals$ and a
constant $A$ such that $d_\Teich(\gamma(h(t)), \gamma'(t)) \le A$ for all
$t \in \reals$. We may assume that $h$ is bilipschitz. It follows that the
relation $\gamma(h(t)) \leftrightarrow\gamma'(t)$ lifts to a
$\pi_1(S)$-equivariant map $X_\gamma\to X_{\gamma'}$, well defined up to
moving points a uniformly bounded distance, and this map is a fiber
respecting quasi-isometry; see Proposition~4.2 of
\cite{FarbMosher:quasiconvex}.

The following theorem was proved independently by Mosher and by
Bowditch:

\begin{theorem}[\cite{Mosher:StableQuasigeodesics}; \cite{Bowditch:stacks}]
\label{TheoremCoboundedGeodesic}
Given a coarse Lipschitz, cobounded, bi-infinite path $\gamma\from \reals
\to\Teich$, the bundle $X_\gamma$ is Gromov hyperbolic if and only if there
exists a cobounded \Teichmuller\ geodesic $\ell$ which fellow travels
$\gamma$.
\end{theorem}

We need a uniform version of Theorem~\ref{TheoremCoboundedGeodesic},
proved in \cite{Mosher:StableQuasigeodesics}:

\begin{proposition}
\label{PropUniformCobounded}
For every cocompact, $\MCG(S)$-equivariant subset $B$ of $\Teich$,
and for every $\delta > 0$, there exist $K \ge 1$, $C \ge 0$ such that if
$\gamma \from \reals \to \Teich$ is a bi-infinite path with image in
$B$, and if $X_\gamma$ is $\delta$ hyperbolic, then there exists a
bi-infinite geodesic $\ell$ in $\Teich$, and a $K$ bilipschitz
homeomorphism $h \from \reals \to \reals$ such that
$d(\gamma(h(t)),\gamma'(t)) \le C$ for all $t \in \reals$. 
\end{proposition}

\paragraph{Singular \solv\ spaces.} Given a cobounded \Teichmuller\ geodesic
$\ell$, there is a natural singular \solv\ metric on $X_\ell$ which we
denote $X^\solv_\ell$. To define the metric, choose a quadratic
differential representing a tangent vector on $\ell$, with horizontal and
vertical measured foliations $(\fol_v,d\mu_v)$, $(\fol_h,d\mu_h)$. On $S
\cross \reals$ we put the singular solv metric $e^{2t} d\mu_v^2 + e^{-2t}
d\mu_h^2 + dt^2$. Lifting to the universal cover we obtain the singular
\solv\ space $X^\solv_\ell$, on which $\pi_1 S$ acts by isometries. Note
that the assignment $\ell \to X^\solv_\ell$ is natural, and so each $f \in
\MCG$ induces an isometry $X^\solv_\ell \to X^\solv_{f(\ell)}$. In
particular, if $\ell$ is the axis of a pseudo-Anosov $g \in \MCG(S)$ then
the extension group $\pi_1(S)\semidirect_g \Z$ acts on $X^\solv_\ell$ by
isometries, properly discontinuously and cocompactly; in this context we
shall also use $X^\solv_g$ to denote $X^\solv_\ell$.

Note also that the identity map $X_\ell \to X^\solv_\ell$ is a
quasi-isometry; see Proposition~4.2 of \cite{FarbMosher:quasiconvex}. When
$\ell$ is an axis this is evident from the fact that the identity map is
equivariant with respect to the action of $\pi_1(S)\semidirect_f \Z$, and
the latter group acts cocompactly by isometries on both $X_\ell$ and
$X^\solv_\ell$.

\paragraph{Orbifolds.} All of our orbifolds will be closed, hyperbolic
\nb{2}orbifolds with only cone points: no mirror edges and hence no
dihedral points. The reason for this restriction is that an orbifold with
mirror edges cannot support a pseudo-Anosov homeomorphism. 

An orbifold $\O$ has a \Teichmuller\ space $\Teich(\O)$ on which the mapping
class group $\MCG(\O) = \Homeo(\O) / \Homeo_0(\O)$ acts, where $\Homeo(\O)$
is the topological group of orbifold homeomorphisms, and $\Homeo_0(\O)$ is
the normal subgroup $\Homeo_0(\O)$ which is the component of the identity
map. 

As with surface mapping class groups, there is a universal extension of
$\pi_1\O$ which can be described by two equivalent short exact sequences,
generalizing the Dehn-Nielsen-Epstein-Baer theorem; we take the description
of this extension from \cite{FarbMosher:quasiconvex}. Unlike surface
mapping class groups, when $\O$ has at least one cone point then the
universal extension group will \emph{not} be isomorphic to the once
punctured mapping class group $\MCG(\O,o)$, where $o \in \O$ is a generic
base point. Instead, the universal extension is the group of lifts of
$\MCG(\O)$ to the universal cover $\wt\O$. To be precise,
$\wt\MCG(\O) = \wt\Homeo(\O) / \wt\Homeo_0(\O)$ where $\wt\Homeo(\O)$ is the
topological group of homeomorphisms of $\wt\O$ that respect orbits of the
action of $\pi_1\O$, and $\wt\Homeo_0(\O)$ is the component of the identity
element of $\wt\Homeo(\O)$ is a normal subgroup; equivalently,
$\wt\Homeo_0(\O)$ consists of those elements of $\wt\Homeo(\O)$ acting
trivially on the circle at infinity of $\wt\O$. As a consequence, we have
an injection $\wt\MCG(\O) \inject \Homeo(S^1\wt\O)$. The analogue
of the Dehn-Nielsen-Epstein-Baer theorem says
$$\xymatrix{
1 \ar[r] & \pi_1(\O) \ar[r] \ar@{=}[d] & \wt\MCG(\O) \ar[r] \ar@2{~}[d]
& \MCG(\O) \ar[r] \ar@2{~}[d] & 1 \\ 
1 \ar[r] & \pi_1(\O) \ar[r] & \Aut(\pi_1(\O)) \ar[r] & \Out(\pi_1(\O))
\ar[r] & 1 }
$$
We say that an automorphism $\hat f \from \Aut(\pi_1(\O,o))$
\emph{represents} a mapping class $f \in \MCG(\O)$ if $\hat f$ maps to
$f$ under the homomorphism $\Aut(\pi_1(\O,o)) \to \Out(\pi_1(\O,o))
\approx \MCG(\O)$.

Note that there is a surjective homomorphism $\MCG(\O,o)
\to\Aut(\pi_1(\O,o))$ defined in the obvious manner. But if $\O$ has any
cone points then this homomorphism has a nontrivial kernel. The kernel is
generated by pushing the base point $o$ around any generic closed curve of
the form $\gamma *\rho * \gamma^\inv$, where $\rho$ is a based closed curve
encircling an order $n$ cone point and going around that point $n$ times,
and $\gamma$ is a path from $o$ to the base point of $\rho$.

There is a canonical hyperbolic surface bundle $X_\Teich$ over the
\Teichmuller\ space $\Teich$ of an orbifold $\O$, on which $\wt\MCG(\O)$
acts by isometries, and Theorem~\ref{TheoremCoboundedGeodesic} and
Proposition~\ref{PropUniformCobounded} are true in this context.

\paragraph{The model space $X_H$.}

Consider an irreducible subgroup $H \subgroup \MCG(S)$. Construct a graph
$G_H$ in \Teichmuller\ space on which $H$ acts properly and cocompactly by
isometries, as follows. Take an orbit of $H$, and use a piecewise geodesic
to connect points which differ by a generator. Since $H$ is finitely
generated and acts properly, we can choose these connecting paths in an
$H$-equivariant way so that they are disjoint except at their endpoints. It
follows that $G_H \subset \Teich$ is an embedded Cayley graph of $H$,
metrized so that each edge is a piecewise geodesic segment in
$\Teich$. Each path in $G_H$ can be regarded as a piecewise geodesic in
$\Teich$. The orbit map from $H$ to $G_H$ is a quasi-isometry, and the
inclusion $G_H \inject\Teich$ is uniformly proper.

By restricting the bundle $X_\Teich \to \Teich$ to $G_H$ we obtain the
canonical $\hyp^2$ bundle $X_H\to G_H$ and a piecewise Riemannian metric
thereon. The group $\Gamma_H$ acts on $X_H$ by isometries, properly
discontinuously and cocompactly, and so any orbit map $\Gamma_H \to X_H$ is
a quasi-isometry.  

Let $\QI_f(X_H)$ denote the subgroup of $\QI(X_H)$ represented by
quasi-isometries $\Phi \from X_H \to X_H$ which uniformly coarsely respect
the fibers, in the sense that there is a constant $A$ such that for each
$x \in G_H$ there exists $x' \in G_H$ such that $d_\Haus(\Phi(D_x),D_{x'})
\le A$. Choosing $x' = \phi(x)$ for each $x \in G_H$, we obtain an induced
quasi-isometry $\phi \from G_H \to G_H$. An orbit map $o \from \Gamma_H \to
X_H$ takes any coset of $\pi_1 S$ in $\Gamma_H$ to with a uniformly finite
Hausdorff distance of some fiber of $X_H$, and every fiber in $X_H$ is has
uniformly finite Hausdorff distance from the image of some coset. It
follows that $o$ induces an isomorphism $\QI_f(\Gamma_H)\approx\QI_f(X_H)$,
and henceforth we will identify these two groups. Note that we obtain a
commutative diagram of homomorphisms
$$\xymatrix{
\QI_f(\Gamma_H) \ar[r]^{\approx} \ar[d] & \QI_f(X_H) \ar[d] \\
\QI(H) \ar[r]^{\approx} & \QI(G_H)
}$$

\paragraph{Hyperbolic lines and coarse axes in $X_H$.}

A bi-infinite path $\gamma$ in $\Teich$ is called a
\emph{hyperbolic line} if $X_\gamma$ is Gromov hyperbolic. By
applying Theorem~\ref{TheoremCoboundedGeodesic}, this is equivalent to the
existence of a cobounded \Teichmuller\ geodesic $\ell$ which fellow travels
$\gamma$ in $\Teich$, and in this case there is an induced fiber respecting
quasi-isometry $X_\gamma \approx X^\solv_\ell$. In this situation, we adopt
the notation $X^\solv_\gamma$ for the metric space $X^\solv_\ell$.

\begin{lemma} 
\label{LemmaCoarseAxis}
If $\gamma$ is a hyperbolic line with image contained in
$G_H$, fellow travelling a cobounded geodesic $\ell \subset \Teich$,
the following are equivalent:
\begin{enumerate}
\item There exists a pseudo-Anosov element $f \in G_H$ such that $\gamma$
is coarsely equivalent in $G_H$ to any orbit of $f$.
\item $\ell$ is the axis of some pseudo-Anosov element $f \in \MCG(S)$.
\end{enumerate}
\end{lemma}

If these happen then we say that $\gamma$ is a \emph{coarse $\Teich$-axis}.

\begin{proof} To prove that (1) implies (2), the axis of $f$ is coarsely
equivalent to any orbit of $f$ and so is coarsely equivalent to $\ell$, but
this implies that the axis of $f$ equals $\ell$, because two cobounded
\Teichmuller\ geodesics which are coarsely equivalent are equal (Lemma~2.4
of \cite{FarbMosher:quasiconvex})

To prove that (2) implies (1), suppose that $\ell$ is the axis of a
pseudo-Anosov element $g \in \MCG(H)$, and so $\ell$ is coarsely contained
in $G_H$. In the word metric on $\MCG(H)$, the infinite cyclic subgroup
$\<g\>$ is coarsely contained in $H$.
Applying Proposition~\ref{PropSubgroups}, some finite index subgroup
$\<g^n\>$ of $\<g\>$ is contained in $H$. It follows that any orbit of
$f=g^n$ in $G_H$ is coarsely equivalent to $\gamma$.
\end{proof}

\begin{theorem}[\cite{FarbMosher:sbf}] 
\label{TheoremSolvEqQIf}
If $\ell,\ell'$ are cobounded geodesics in $\Teich$, if $\ell$ is the axis
of a pseudo-Anosov mapping class, and if $\Phi\from X_\ell \to X_{\ell'}$
is a fiber respecting quasi-isometry, then $\ell'$ is the axis of a
pseudo-Anosov mapping class and $\Phi$ is a bounded distance from an
isometry $X^\solv_\ell \to X^\solv_{\ell'}$.  
\end{theorem}

\begin{corollary} 
\label{CorollarySolvIsometry}
If $\gamma$ is a hyperbolic line in $G_H$ and $\Phi \in
\QI_f(X_H)$ then $\gamma = \Phi \composed \gamma$ is a hyperbolic line in
$G_H$. Moreover, $\gamma$ is a coarse axis if and only if $\gamma'$ is a
coarse axis, in which case the fiber respecting quasi-isometry $X_\gamma
\to X_{\gamma'}$ induced by $\Phi$ is a bounded distance from an isometry
$X^\solv_\gamma \to X^\solv_{\gamma'}$.
\end{corollary}

\begin{proof} Suppose $\gamma$ is a hyperbolic line in $G_H$, that is,
$X_\gamma$ is a Gromov hyperbolic metric space. The quasi-isometry $\Phi
\from X_H \to X_H$ restricts to a quasi-isometry $X_\gamma \to
X_{\gamma'}$, and so $X_{\gamma'}$ is Gromov hyperbolic, that is, $\gamma'$
is a hyperbolic line.

By using Theorem~\ref{TheoremCoboundedGeodesic} we may replace $X_\gamma$
and $X_{\gamma'}$ with $X^\solv_\ell$ and $X^\solv_{\ell'}$ for
\Teichmuller\ geodesics $\ell,\ell'$ fellow travelling $\gamma,\gamma'$,
respectively, and then we apply Theorem~\ref{TheoremSolvEqQIf}. 
\end{proof}

\section{The group $\QI_f(\Gamma_H) \approx \QI_f(X_H)$}

\paragraph{The quasi-symmetry group $\QSym(S^1)$.} 

The quasi-isometry group of the hyperbolic plane $\hyp^2$ acts
faithfully on the boundary circle $S^1$, and the image of this
action is precisely the group of quasi-symmetric homeomorphisms of
$S^1$, denoted $\QSym(S^1)$. 

Consider $X_H$. For each $x \in G_H$, the fiber over $x$ is $D_x \subset
X_H$. Fix a base point $p_0 \in G_H$ and let $D_0
= D_{p_0}$. Fix  an isometric identification of $D_0$ with $\hyp^2$, so the
Gromov boundary of
$D_0$ is identified with $S^1$ and $\QI(D_0) \approx \QSym(S^1)$. 

There is an induced homomorphism $\theta \from \QI_f(X_H)\to
\QSym(S^1)$ defined as follows. For each $\Phi \in \QI_f(X_H)$, the
set $\Phi(D_0)$ is Hausdorff close to $D_0$ and so a closest point map
$\Phi(D_0) \to D_0$, precomposed with $\Phi$, induces a quasi-isometry $D_0
\to D_0$ whose class $\theta(\Phi) \in \QI(\hyp^2) = \QSym(S^1)$ is
well-defined independent of the choice of a closest point map.

\begin{claim} The homomorphism $\theta \from \QI_f(X_H) \to\QSym(S^1)$ is
injective.
\end{claim}

\begin{proof} The proof is in principle the same as in
\cite{FarbMosher:sbf}, although there it was couched in terms of the limit
set of a Schottky group. We shall reformulate the proof without
referring to Schottky groups, as follows.

Consider a pseudo-Anosov element $g \in H$ with coarse axis $\gamma$, and
choose a fiber $D_g \subset X^\solv_g$. Since $D_0$ and $D_g$ have finite
Hausdorff distance, there is a canonical homeomorphism $\bdy D_g \approx
\bdy D_0 = S^1$. With respect to this homeomorphism, let $E_s(g) \subset
S^1$ be the set of endpoints of leaves of the stable foliation on $D_g$,
let $E_u(g) \subset S^1$ be the set of endpoints of unstable leaves, and let
$E(g) = E_s(g) \union E_u(g)$, a disjoint union.  Note that if $g,g' \in H$
are pseudo-Anosov then either $\VN(g)=\VN(g')$ and $E(g)=E(g')$, or $\VN(g)
\intersect \VN(g') = \{\Id\}$ and $E(g)\intersect E(g') = \emptyset$. 

Fix a conjugacy class $\C$ of pseudo-Anosov elements of $H$, fix a coarse
axis $\gamma_g$ for one particular element of $H$, and for any $h \in H$ we
obtain a coarse axis $h(\gamma_g)$ for $hgh^\inv$. There is a fixed $\delta
> 0$ so that the spaces $X^\solv_g$ are all $\delta$-hyperbolic for $g \in
\C$. Applying Proposition~\ref{PropUniformCobounded} with $B = G_H$, it
follows that the axis $\ell_g \subset \Teich$ for $g$ is uniformly Hausdorff
close to $\gamma_g$. 

Let $\Phi \from X_H \to X_H$ be a fiber respecting quasi-isometry such that
$\theta(\Phi)$ is the identity. Let $\phi \from G_H \to G_H$ be the
quasi-isometry induced by $\Phi$. Applying
Proposition~\ref{CorollarySolvIsometry}, for each $g
\in\C$ the image $\phi(\gamma_g)$ is a coarse axis for some pseudo-Anosov
element $g' \in H$, and the induced map $X_{\gamma_g} \to X_{\gamma_{g'}}$
is a quasi-isometry with uniform constants. It follows that
$X_{\gamma_{g'}}$ is $\delta'$ hyperbolic for $\delta'$ independent of $g$.
Applying Proposition~\ref{PropUniformCobounded}, $\gamma_{g'}$ is a
uniformly finite Hausdorff distance from the axis $\ell_{g'}$ of $g'$.
We have canonical identifications $\bdy D_g \approx \bdy D_{g'} \approx \bdy
D_0 = S^1$, and the map $\bdy D_g \to\bdy D_{g'}$ induced by $\Phi$ must
agree with the identity map on $\bdy D_0$, since $\theta(\Phi)$ is the
identity. On the other hand, $\Phi(E(g)) = E(g')$, implying that $E(g) =
E(g')$, and so $\VN(g) = \VN(g')$. It follows that $\ell_g = \ell_{g'}$, and
so $\phi(\gamma)$ is coarsely equivalent to $\gamma$, with a uniformly
finite Hausdorff distance. 

The collection of axes $\gamma_g$ for $g \in \C_x$ coarsely separate points
in $G_H$. That is, there exists constants $r,s \ge 0$ so that for each $x\in
G_H$ there is a finite subset $\C_x \subset \C$ such that $x$ is within
distance $r$ of each coarse axis $\gamma_g$ for $g \in \C_x$, and the set
of points within distance $r$ of these coarse axes has diameter at most
$s$. Since $\phi(\gamma_g)$ is uniformly coarsely equivalent to $\gamma_g$
for each $g\in \C$, it follows that $d(\phi(x),x)$ is uniformly finite, in
other words the quasi-isometry $\phi \from G_H \to G_H$ is equivalent to the
identity. 

For each $x \in G_H$, it follows that $\Phi(D_x)$ is a uniformly
finite distance from $D_x$, and since the boundary values of the closest
point projection $\Phi(D_x) \to D_x$ must agree with the identity map on
$S^1$ it follows that the closest point projection is a uniformly finite
distance from the identity on $D_x$. This proves that $\Phi$ is equivalent
to the identity map on $X_H$.
\end{proof}

Having proved that $\theta$ is injective, we will identify $\QI_f(X_H)
\approx \QI_f(\Gamma_H)$ with its image in $\QSym(S^1)$.

There are several other subgroups of $\QSym(S^1)$ that are of
significance to us. In particular, for each pseudo-Anosov $f \in H$ with
coarse axis $\gamma$ fellow travelling a \Teichmuller\ axis $\ell$ for $f$,
we have constructed a fiber preserving quasi-isometry $X_\gamma \approx
X^\solv_\gamma = X^\solv_\ell$. The group $\Isom(X^\solv_\gamma)$ takes
fibers to fibers. If we choose a fiber $D_\gamma \subset
X^\solv_\gamma$ then each isometry of $X^\solv_\gamma$ takes $D_\gamma$ to
another fiber coarsely equivalent to $D_\gamma$, and composing with the
closest point projection go $D_\gamma$ we obtain a self quasi-isometry of
$D_\gamma$, thereby constructing an injection $\Isom(X^\solv_\gamma) \to
\QI(D_\gamma)$. However, $D_0$ and $D_\gamma$ are coarsely equivalent in
$X_H$ and so we obtain an injection $\Isom(X^\solv_\gamma) \inject
\QSym(S^1)$; again we identify $\Isom(X^\solv_\gamma)$ with the image of
this injection.

\paragraph{Mapping classes which descend to subcovers.} Consider an
orbifold covering map $\pi \from \Qorb \to \Porb$, and fix generic base
points $q \in \Qorb$, $p=\pi q \in \Porb$, inducing an injection of
orbifold fundamental groups $\pi_1(\Qorb,q) \to
\pi_1(\Porb,p)$. A mapping class $f \in \MCG(\Qorb)$ \emph{descends} to a
mapping class $g \in \MCG(\Porb)$ if there exist orbifold homeomorphisms
$F\from \Qorb \to\Qorb$ and $G \from\Porb \to\Porb$ representing $f$ and
$g$, respectively, such that $F$ is a lift of $G$, that is, $F\composed\pi
=\pi\composed G$. 

\begin{lemma}
\label{LemmaDescent}
Consider $\pi \from \Qorb \to \Porb$, and $q \in \Qorb$, $\pi q = p \in
\Porb$ as above. Given $f \in \MCG(\Qorb)$, $g\in \MCG(\Porb)$, the
following are equivalent:
\begin{enumerate}
\item $f$ descends to $g$. 
\item There exist automorphisms $\hat f \from \pi_1(\Qorb,q) \to
\pi_1(\Qorb,q)$ and $\hat g \from \pi_1(\Porb,p) \to \pi_1(\Porb,p)$
representing $f,g$, respectively, such that $\hat f$ extends to $\hat g$
(with respect to the injection $\pi_1(\Qorb,q) \inject \pi_1(\Porb,p)$). 
\item For any automorphism $\hat f  \from \pi_1(\Qorb,q) \to
\pi_1(\Qorb,q)$ representing $f$, there exists an automorphism $\hat g 
\from \pi_1(\Porb,p) \to\pi_1(\Porb,p)$ representing $g$, such that
$\hat f$ extends to $\hat g$.
\end{enumerate}
Moreover, if $f \in \MCG(\Qorb)$ descends to $\MCG(\Porb)$, then the
element $g \in \MCG(\Porb)$ to which it descends is unique.
\end{lemma}

\begin{proof} In this proof we will use the fact that if $\O$ is an
orbifold with base point $x$ then we have a surjective homomorphism
$\MCG(\O,x) \to \wt\MCG(\O,x) \approx \Aut(\O,x)$, as described earlier.

Obviously (3) implies (2). 

To prove (2) implies (1), suppose that $\hat f, \hat g$ exist as in (2).
Choose an orbifold homeomorphism $G \from (\Porb,p) \to (\Porb,p)$ that
induces $\hat g$. By the lifting lemma in the category of orbifolds, we may
lift $G$ to a homeomorphism $F \from (\Qorb,q) \to (\Qorb,q)$ that
induces~$\hat f$. Forgetting base points, clearly $G$ represents $g$ and
$F$ represents $f$.

To prove (1) implies (3), suppose there exist
orbifold homeomorphisms $F\from\Qorb \to \Qorb$ and $G \from\Porb \to \Porb$
representing $f$ and $g$, respectively, such that $F$ is a lift of $G$. We
may isotop $G$ so that it fixes $p$, and we may lift this to an isotopy of
$F$ so that $F$ permutes $\pi^\inv p$, but it may not be true that $F(q)=q$.
If not, pick a path $\gamma$ in $\Qorb$ from $F(q)$ to $q$ whose projection
$\pi\gamma$ passes only through generic points of $\Porb$, and postcompose
$G$ by the isotopy of $\Porb$ that pushes $p$ around $\pi \gamma$ to obtain
a new $G$. This lifts to an isotopy of $\Qorb$ that pushes $F(q)$ along
$\gamma$ to $q$, and postcomposing $F$ we obtain a new $F$ so that
$F(q)=q$. It follows that $F$ and $G$ induce $\hat f \in
\Aut(\pi_1(\Qorb,q))$ and $\hat g \in \Aut(\pi_1(\Porb,p))$ satisfying the
requirements of the lemma.

To prove uniqueness of $g$ given $f$, for each automorphism $\hat f \from
\pi_1(\Qorb,q) \to\pi_1(\Qorb,q)$ representing $f$ it suffices to prove
that there is at most one automorphism $\hat g \from \pi_1(\Porb,p)
\to\pi_1(\Porb,p)$ to which $\hat f$ extends. We have an isomorphism
$\Aut(\pi_1(\Porb,p)) \approx \wt\MCG(\Porb)$, and an injection
$\wt\MCG(\Porb) \inject \QSym(S^1)$; in other words, an
automorphism of $\pi_1(\Porb,p)$ is determined by its boundary values. But
the finite index injection $\pi_1(\Qorb,q) \inject \pi_1(\Porb,p)$ induces
a homeomorphism of boundaries, and so the boundary values of an extension
of $\hat f$ are completely determined by $\hat f$.
\end{proof}

Given an orbifold covering map $\Qorb \to \Porb$, a subgroup $A \subset
\MCG(\Qorb)$ \emph{descends} to a subgroup $B \subset \MCG(\Porb)$ if there
is an isomorphism $A \approx B$ such that each element of $A$ descends to
the corresponding element of $B$. We also say simply that $A$ descends to
$\Porb$.

In \cite{FarbMosher:sbf}, given a surface $S$ and a free subgroup $H
\subgroup \MCG(S)$, we studied the problem of finding a smallest orbifold
subcovering map $S \to \O$ such that $H$ descends to $\O$. Since $H$
was free, it was more or less obvious that $H$ descends if and only if the
elements of some free basis descend, and this allowed us to ignore some
subtleties of the concept of ``descent''. In the present situation we are
working with arbitrary subgroups, and the question arises whether one can
check descent on an element by element basis. The answer turns out to be
yes, as we show in Lemma~\ref{LemmaExtension} below; the proof is based on
the uniqueness of descent for mapping classes proved of
Lemma~\ref{LemmaDescent}. Uniqueness of descent for mapping classes also
implies uniqueness of descent for subgroups, as we now see.

\begin{lemma}
\label{LemmaExtension}
Consider $\pi \from \Qorb \to \Porb$, $q \in \Qorb$, $p=\pi q \in \Porb$ as
above. For any subgroup $A \subgroup \MCG(\Qorb)$, there is at most one
subgroup of $\MCG(\Porb)$ to which $A$ descends. Moreover, the following
are equivalent:
\begin{enumerate}
\item $A$ descends to some subgroup of $\MCG(\Porb)$.
\item Each element of $A$ descends to some element of $\MCG(\Porb)$.
\item There exists a generating set of $A$ such that each generator
descends to some element of $\MCG(\Porb)$.
\end{enumerate}
\end{lemma}

\begin{proof} Obviously (1) implies (2) implies (3). To prove (3) implies
(1), let $G \subset A$ be a generating set such that each $g \in G$
descends to $\MCG(\Porb)$, equivalently, we can choose a representative
$\hat g \in \Aut(\pi_1(\Qorb,q))$ that extends to some $\hat f \in
\Aut(\pi_1(\Porb,p))$. The composition of each word in the $\hat g$
therefore extends to some element of $\Aut(\pi_1(\Porb,p))$. We must prove
that if $w(g)$ is a word in the letters $g \in G$ that represents the
identity element of $G$, if $w(\hat g)$ is the corresponding word in the
letters $\hat g$, and if $w(\hat f)$ is the corresponding word in the
letters $\hat f$, then $w(\hat f)$ is a representative of the identity
element of $\MCG(\Porb)$. But we know that $w(\hat g)$ is a representative
of the identity element of $\MCG(\Qorb)$, which means that $w(\hat g)$
agrees with the inner automorphism of $\pi_1(\Qorb,q)$ defined by some
element $\gamma \in \pi_1(\Qorb,q)$. The word $w(\hat f)$ represents an
extension of $w(\hat g)$ to an automorphism $\pi_1(\Porb,p)$, but $\gamma$
acts by conjugation on $\pi_1(\Porb,p)$ since $\pi_1(\Qorb,q) \subgroup
\pi_1(\Porb,p)$, and the actions of $w(\hat f)$ and $\gamma$ agree on the
circle at infinity because they agree with restricted to $\pi_1(\Qorb,q)$.
It follows that $w(\hat f)$ and the conjugation action of $\gamma$ are the
same automorphism of $\pi_1(\Porb,p)$. This proves that $w(\hat f)$ is an
inner automorphism of $\pi_1(\Porb,p)$, and so $w(\hat f)$ is a
representative of the identity element of $\MCG(\Porb)$.
\end{proof}

Given a pseudo-Anosov $f \in \MCG(S)$, in \cite{FarbMosher:sbf} we defined
the ``maximal orbifold subcover'' $S \to \O_f$ to which $f$ descends. The
orbifold $\O_f$ is constructed explicitly as follows. Let
$\ell \subset
\Teich$ be the axis of $f$. Let $X^\solv_\ell$ be the singular
\solv-manifold associated to $\ell$. Pick a fiber $D_\ell$, a singular
Euclidean space, with a horizontal and vertical measured foliation, on
which $\pi_1 S$ acts properly and cocompactly by isometries preserving the
two foliations. Consider the full group of isometries
$\Isom(X^\solv_\ell)$. The fibers of $X^\solv_\ell$ are preserved by
$\Isom(X^\solv_\ell)$, and so the action of $\Isom(X^\solv_\ell)$ descends
to an action on $\ell$ which is discrete and cocompact. The quotient group
acting on $\ell$ is denoted $C_\ell$, and it is isomorphic to either the
infinite cyclic group $\Z$ or the infinite dihedral group
$D_\infinity$. The kernel of this action, denoted $\Isom_f(X^\solv_\ell)$,
the subgroup of $\Isom(X^\solv_\ell)$ that preserves each fiber, and the
restriction of $\Isom_f(X^\solv_\ell)$ to $D_\ell$ is the full group of
isometries of $D_\ell$ that preserve the horizontal and vertical measured
foliations. The group $\Isom_f(X^\solv_\ell)$ acts properly and cocompactly
on $D_w$, containing $\pi_1(S)$ as a finite index subgroup. It follows that
$\O_f = \O_\ell = D_\ell / K_\ell$ is an orbifold subcover of $S$; let $\pi
\from S \to \O_f$ be the projection. By construction, for any pseudo-Anosov
homeomorphism $F$ representing $f$ there is a pseudo-Anosov homeomorphism
$G$ of $\O_f$ such that $\pi\composed F = G \composed \pi$, and so $f$
descends to $\O_f$. 

To explain maximality of $\O_f$, suppose that $f$ descends through a
subcover $S\to \Qorb$ to some $g \in \MCG(\Qorb)$. Choose $G\from \Qorb \to
\Qorb$ representing $g$ and a lift $F\from S \to S$ representing $f$. Since
$f$ is irreducible, $g$ must also be irreducible, and so $g$ has a
pseudo-Anosov representative $G'$. The isotopy from $G$ to $G'$ lifts to an
isotopy from $F$ to a homeomorphism $F' \from S \to S$ which is obviously
pseudo-Anosov. It follows that $X^\solv_f$ and $X^\solv_g$ are isometric. We
obtain inclusions $\pi_1(S) \subgroup \pi_1(\Qorb)
\subgroup \Isom_f(X^\solv_f)$, and so the orbifold covering map $S \to \O_f$
factors as orbifold covering maps $S \to \Qorb \to \O_f$, which gives
maximality of $\O_f$.

\begin{proposition}
If $H\subgroup\MCG(S)$ is irreducible then there exists an orbifold
covering map $\pi \from S \to \O_H$ to which $H$ descends, so that $\O_H$ is
maximal: any orbifold covering map $S \to\Qorb$ to which
$H$ descends composes with a covering map $\Qorb \to\O_H$ to give $\pi$.
The group $\pi_1\O_H$ can be characterized in
$\QSym(S^1)$ as
$$\pi_1 \O_H = \bigcap_g \pi_1\O_g, \quad\text{over all pseudo-Anosov $g
\in H$.}
$$
\end{proposition}

\begin{proof}
In the elementary case, where $H \subgroup \VN(f)$ for some pseudo-Anosov
$f \in \MCG(S)$, we clearly have $\O_H = \O_f$.

In the nonelementary case, $H$ is generated by its pseudo-Anosov elements.
Each pseudo-Anosov $f \in H$ descends through a maximal subcover $S \to
\O_f$. We have inclusions $\pi_1 S \subgroup \pi_1 \O_f \subgroup
\QSym(S^1)$, the first of which has finite index. Define
$\pi_1\O_H$ to be the intersection of $\pi_1\O_f$ over all pseudo-Anosov $f
\in H$. Define $\O_H$ to be the quotient orbifold of $\pi_1\O_H$.
The action of $H$ by conjugation on $\QSym(S^1)$ permutes the subgroups
$\pi_1\O_f$, for pseudo-Anosov $f \in H$. In particular the action of
any pseudo-Anosov $g \in H$ permutes the $\pi_1\O_f$, and so $g$ preserves
$\pi_1\O_H$, that is, $g$ descends to $\O_H$. Applying
Lemma~\ref{LemmaExtension}, $H$ descends to $\O_H$. 

To prove maximality, consider an orbifold covering $S \to \Qorb$ such that
$H$ descends to $\Qorb$, and so each pseudo-Anosov $f \in H$ descends
to $\Qorb$, which may therefore be composed with a covering map $\Qorb \to
\O_f$ to obtain the covering map $S \to \O_f$. It follows that $\pi_1(S)
\subgroup \pi_1(\Qorb) \subgroup \pi_1\O_f$ in $\QSym(S^1)$, for
all pseudo-Anosov $f \in H$, implying that $\pi_1(\Qorb) \subgroup
\pi_1(\O_H)$, which implies that $S \to \Qorb$ composes with a covering map
$\Qorb \to \O_H$ to yield $S \to \O_H$.
\end{proof}

\paragraph{Computation of $\QI_f(X_H)$.} Consider an irreducible subgroup
$H \subgroup \MCG(S)$. Let $\O=\O_H$. The subgroup $H
\subgroup \MCG(S)$ descends to a subgroup of $\MCG(\O)$ which will
be denoted $H'$. Let $\C$ denote the relative commensurator of $H'$ in
$\MCG(\O)$. We have a commutative diagram of extensions
$$\xymatrix{
1 \ar[r] & \pi_1(\O) \ar[r] \ar@{=}[d] & \Gamma_{H'} \ar[r] \ar[d] & H'
\ar[r] \ar[d] & 1 \\ 
1 \ar[r] & \pi_1(\O) \ar[r] \ar@{=}[d] & \Gamma_\C \ar[r] \ar[d] & \C
\ar[r] \ar[d] & 1 \\ 
1 \ar[r] & \pi_1(\O) \ar[r] & \wt\MCG(\O) \ar[r] & \MCG(\O) \ar[r]
& 1 
}
$$
in which all vertical arrows are inclusions. We also have a commutative
diagram
$$\xymatrix{
1 \ar[r] & \pi_1(\O) \ar[r] \ar[d] & \Gamma_{H'} \ar[r] \ar[d] & H'
\ar[r] \ar@{=}[d] & 1 \\ 
1 \ar[r] & \pi_1(S) \ar[r] & \Gamma_H \ar[r] & H \ar[r]
& 1 
}
$$
where the vertical arrows are inclusions with finite index image. We may
therefore regard $\Gamma_{H'}$ as a finite index subgroup of
$\Gamma_H$, and so the inclusion induces a quasi-isometry
$\QI_f(\Gamma_{H'}) \approx \QI_f(\Gamma_H)$. 

Since $\C$ is the relative commensurator of $H'$ in $\MCG(\O)$, it follows
that $\Gamma_\C$ is the relative commensurator of $\Gamma_{H'}$ in
$\wt\MCG(\O)$, and so we obtain a homomorphism $\Gamma_\C \to
\QI(\Gamma_{H'})$; to be precise, conjugation by $g \in \Gamma_\C$ maps some
finite index subgroup of $\Gamma_{H'}$ isomorphically to another finite
index subgroup, inducing a quasi-isometry of $\Gamma_{H'}$. Moreover,
clearly $\Gamma_\C$ preserves cosets of $\pi_1\O$, and so we actually
obtain a homomorphism $\Gamma_\C\to \QI_f(\Gamma_{H'})\approx
\QI_f(\Gamma_H)$. 

To see that this homomorphism is injective, by postcomposing with
$\QI_f(\Gamma_H) \to \QSym(S^1)$ we obtain a homomorphism
$\Gamma_\C \to \QSym(S^1)$ which factors through two injections $\Gamma_\C
\to \wt\MCG(\O) \to \QSym(S^1)$.

Now we prove that the homomorphism $\Gamma_\C \to \QI_f(\Gamma_H)$ is
surjective: for each fiber respecting quasi-isometry $\Phi \from X_H \to
X_H$ there exists $F \in\Gamma_\C$ such that $\Phi$ agrees within
uniformly bounded distance with the action of $F$. Let $\phi \from G_H \to
G_H$ be the quasi-isometry induced by $\Phi$.

Consider a pseudo-Anosov $g \in H$ with coarse axis $\gamma \subset G_H$.
The image $\phi(\gamma)$ is a coarse axis for a pseudo-Anosov $g'
\in H$. $\Phi$ induces a fiber respecting quasi-isometry $\Phi
\from X_\gamma \to X_{\phi(\gamma)}$ which, by
Corollary~\ref{CorollarySolvIsometry}, is a bounded distance from an
isometry $X^\solv_g \to X^\solv_{g'}$, and so
$\Phi$ conjugates $\Isom(X^\solv_g)$ to $\Isom(X^\solv_{g'})$. Restricting
to any fiber, $\Phi$ conjugates $\Isom_f(X^\solv_g) =
\pi_1(\O_g)$ to $\Isom_f(X^\solv_{g'}) = \pi_1(\O_{g'})$. We may regard
this conjugation as taking place in $\QSym(S^1)$.

We have proved that the conjugation action of $\Phi$ in
$\QSym(S^1)$ permutes the collection of groups $\pi_1 \O_g$, over
all pseudo-Anosov elements $g \in H$. It follows that $\Phi$ preserves the
subgroup
$$\pi_1 \O_H = \bigcap_g \pi_1\O_g
$$
that is, $\Phi \in \Aut(\pi_1\O_H) \approx \wt\MCG(\O_H)$. We shall use $F$
to denote $\Phi$ regarded as an element of the group $\wt\MCG(\O_H)$, and
$f$ to be its image in the quotient $\MCG(\O_H)$.

Having identified $\Phi$ as an element $F \in \wt\MCG(\O_H)$, we must now
show that $F \in \Gamma_{\C}$, that is, $F \Gamma_{H'} F^\inv$ and
$H'$ are commensurable in $\wt\MCG(\O_H)$. Passing to the quotient group
$\MCG(\O_H)$, it suffices to prove that $f H' f^\inv$ and $H'$ are
commensurable in $\MCG(\O_H)$. Applying Proposition~\ref{PropSubgroups},
we are reduced to showing that $f H' f^\inv$ is coarsely
equivalent to $H'$ in $\MCG(\O_H)$. Because we have a uniformly proper
embedding $\MCG(\O_H) \to \Teich(\O_H)$ given by any orbit map, and because
$G_{H'}$ is coarsely equivalent to any orbit of the action of $H'$ on
$\Teich(\O_H)$, it suffices to prove that $f(G_{H'})$ is coarsely
equivalent to $G_{H'}$ in $\Teich(\O_H)$.

We regard $\Phi$ as a quasi-isometry of $X_{H'}$. Let $P$ denote a fixed
pseudo-Anosov conjugacy class in $H'$. For each $g \in P$ pick a coarse
axis $\gamma_g \subset G_{H'}$ in an $H'$-equivariant manner, and let $L =
\{\gamma_g \suchthat g \in P\}$. The metric spaces $X_{\gamma_g}$ are all
isometric for $g \in P$; let $\delta$ be a hyperbolicity constant for each
of them. The action of $H'$ permutes the collection of coarse axes
$L$, and by cocompactness of the $H'$ action on $G_{H'}$ it follows that
the union of the coarse axes in $L$ is coarsely equivalent to $G_{H'}$.

The map $\Phi$ induces a fiber respecting quasi-isometry from
$X_{\gamma_g}$ to $X_{\phi(\gamma_g)}$, with uniform quasi-isometry
constants, and taking each fiber to within a uniform Hausdorff distance of
some other fiber. It follows that there exists $\delta'>0$ so that for each
$g\in P$, the space $X_{\phi(\gamma_g)}$ is $\delta'$-hyperbolic.
Taking $B=G_{H'}$ and applying Proposition~\ref{PropUniformCobounded},
there is a constant $D$ such that $\phi(\gamma_g)$ is within Hausdorff
distance $D$ of the axis $\ell_{g'}$ of some pseudo-Anosov element $g' \in
\MCG(\O_H)$. Similarly, replacing $\delta'$ by $\max\{\delta,\delta'\}$ it
follows that $\gamma_g$ is within Hausdorff distance $D$ of the axis
$\ell_g$.

For $g \in P$, the mapping class $f \in \MCG(\O_H)$ takes $\ell_g$ to
some axis in $\Teich(\O_H)$, and from the construction of $f$ the only
candidate is the axis $\ell_{g'}$. It follows that $f(\gamma_g)$ has a
uniformly finite Hausdorff distance from $\phi(\gamma_g)$. But the union of
the coarse axes $\gamma_g \in L$ is coarsely equivalent to $G_{H'}$, as is
the union of the coarse axes $\phi(\gamma_g)$. It follows that $f(G_{H'})$
is coarsely equivalent to $G_{H'}$, as required.

\paragraph{The quantitative version of Corollary~\ref{CorollaryQIGroup}.} 
Recall from the introduction that we need the following
quantitative version of Corollary~\ref{CorollaryQIGroup}: for all $K \ge
1$, $C \ge 0$ there exists $A \ge 0$ such that if $\phi \from \MCG(S,p) \to
\MCG(S,p)$ is a $K,C$ quasi-isometry then there exists $F \in \MCG(S,p)$
such that $d(\phi(G),FG) \le A$ for all $G \in \MCG(S,p)$. The proof of
Corollary~\ref{CorollaryQIRigidity} then follows by a standard argument
first found in \cite{Schwartz:RankOne}; see also
\cite{Mosher:Durham03}. The quantitative version of
Corollary~\ref{CorollaryQIGroup} follows from quantitative versions of
Whyte's Theorem~\ref{TheoremWhyte} and Theorem~\ref{TheoremMCG}.

For Whyte's Theorem~\ref{TheoremWhyte}, the quantitative version says that
for all $K\ge 1$, $C \ge 0$ there exists $R \ge 0$ such that if $\phi \from
\MCG(S,p)\to\MCG(S,p)$ is a $K,C$ quasi-isometry then each coset of $\pi_1
S$ is taken by $\phi$ to within Hausdorff distance $R$ of another coset of
$\pi_1 S$; in this case $\phi$ is said to \emph{$R$-coarsely respect} the
cosets of $\pi_1 S$. For the proof, see \cite{Mosher:Durham03}.

For Theorem~\ref{TheoremMCG} the quantitative version says that for all $K
\ge 1$, $C\ge 0$, $R \ge 0$ there exists $A \ge 0$ such that if $\phi \from
\MCG(S,p)\to\MCG(S,p)$ is a $K,C$ quasi-isometry that $R$-coarsely respects
the cosets of $\pi_1 S$ then there exists $F \in \MCG(S,p)$ such that
$d(\phi(G),FG) \le A$ for all $G \in \MCG(S,p)$.

Consider the argument given above for the computation of $\QI_f(\Gamma_H)$, 
specialized to the case where $H=\MCG(S)$, and so $\O_H=S$. We continue to
use the notation $G_H$ for a piecewise geodesic Cayley graph equivariantly
embedded in $\Teich$, and $X_H$ for the canonical $\hyp^2$ bundle over
$G_H$. In the coarse of this argument we produced a mapping class $F \in
\MCG(S,p)$, with quotient $f \in \MCG(S)$, which is characterized by the
following property: for each pseudo-Anosov $g \in \MCG(S)$, with coarse
axis $\gamma_g \subset G_H$ and axis $\ell_g \subset \Teich$, letting $g' =
fgf^\inv$ with axis $\ell_{g'} = f(\ell_g) \subset \Teich$, we may conclude
that $\gamma_g$ and $\phi(\gamma_g)$ are a bounded distance in $\Teich$ from
$\ell_g$ and $\ell_{g'}$, respectively, and the quasi-isometry
$X_{\gamma_g} \to X_{\phi(\gamma_g)}$ induced by $\Phi$ is a bounded
distance from the isometry $X^\solv_{\ell_g} \to X^\solv_{\ell_{g'}}$
induced by $F$. The bounds in this conclusion depend on the quasi-isometry
constants $K,C$ for $\Phi$ as well as on the constant $R$ for coarse
preservation of fibers, but the bounds also depend on the pseudo-Anosov
$g$ itself, because the proof invoked Proposition~\ref{PropUniformCobounded}
which depends on the hyperbolicity constant for $X_{\gamma_g}$. 

However, we can get around this dependence by the same trick used already
several times: fix a pseudo-Anosov conjugacy class $P$, and choose a coarse
axis $\gamma_g$ for each $g \in P$ in an $\MCG(S)$-equivariant manner, and
it follows that the metric spaces $X_{\gamma_g}$ are all isometric to each
other and hence have the same hyperbolicity constant $\delta$. This
constant $\delta$ thus depends only on the genus of the surface $S$, and
hence the constants that come in the conclusion of
Proposition~\ref{PropUniformCobounded} depend only on $K,C,R$, and the
genus of $S$.

Using the additional fact that each point of $G_H$ is a uniformly bounded
distance from the axis of some $g \in P$, by a bound which again depends
only on the genus of $S$, we obtain the desired conclusion that
$d(\phi(G),FG) \le A$ for all $G \in \MCG(S,p)$, where $A$ depends only on
$K,C,R$, and the genus of $S$.


\providecommand{\bysame}{\leavevmode\hbox to3em{\hrulefill}\thinspace}

\noindent
Lee Mosher:\\
Department of Mathematics, Rutgers University, Newark\\
Newark, NJ 07102\\
E-mail: mosher@andromeda.rutgers.edu

\end{document}